\magnification=1200
\hfuzz=1pt
\def\sqr#1#2{{\vcenter{\vbox{\hrule height.#2pt 
             \hbox{\vrule width.#2pt height#1pt \kern#1pt \vrule width.#2pt} 
             \hrule height.#2pt}}}}

\def\qed{\hfill{\vbox{\hrule\hbox{\vrule\kern3pt
                \vbox{\kern6pt}\kern3pt\vrule}\hrule}}}

\def\ll{{\lambda}}
\def\kk{{\kappa}} 
\def\phi{{\varphi}}

\def\aa{{\alpha}}
\def\bb{{\beta}}
\def\gg{{\gamma}}
\def\GG{{\Gamma}}
\def\Q{{\bf Q}}
\def\R{{\bf R}}
\def\I{{\cal I}}
\def\P{\noindent{\bf Proof.} }
\noindent P. Komj\'ath, Dept. Comp. Sci. E\"otv\"os University,
Budapest, M\'uzeum krt 6--8, 1088, Hungary, {\tt e-mail: kope@cs.elte.hu}

\noindent S. Shelah, Inst. of Mathematics, Hebrew University,
Jerusalem, Israel, 

\noindent{\tt e-mail: shelah@sunrise.huji.ac.il}
\bigskip
\centerline{\bf On uniformly antisymmetric functions}
\medskip
\centerline{\bf 0. Introduction}
\medskip
\noindent Recently there has been\vfootnote{}{No. 502 on
the second author's list. Supported by the Hungarian OTKA grant
No. 1908 and by the grant of the Israeli
Academy of Sciences.\hfil\break\indent AMS subject classification
(1991): 26 A 15, 03 E 50, 04 A 20.} considerable research on symmetric
properties of functions, i.e., when e.g.~continuity is replaced
by the limit properties of $f(x+h)-f(x-h)$ ($h\to 0$). The
excellent monograph [6] surveys most of the recent developments.

The following definiton was considered    by Evans and
Larson (in Santa Barbara, 1984) and by Kostyrko (in Smolenice,
1991).
\medskip
\noindent {\bf Definition.} A {\sl uniformly antisymmetric function}
is an
$f:\R\to\R$ such that for every $x\in\R$ there is a
$d(x)>0$ so that $0<h<d(x)$ implies $\vert f(x+h)-f(x-h)\vert
\geq d(x)$. 
\medskip
They posed the question if  there exists a 
uniformly antisymmetric function.
Kostyrko showed that no such function with a two
element range exists, that is, there is no uniformly
antisymmetric {\sl set} (see [5]). This was extended to
functions with 3-element ranges by Ciesielski in [1]. In [2] 
a uniformly antisymmetric
function $f:\R\to \omega$ was constructed. It had the stronger
property that for every $x\in \R$ the set
$S_x=\{h>0:f(x-h)=f(x+h)\}$ is finite. [2] contains several
other relevant results and questions. Kostyrko's result is
extended to functions defined on any uncountable subfield of the
reals. The authors of [2] ask if this can be extended to countable
subfields, as well. As for functions defined on $\R$ they ask if
there is an $f:\R\to\omega$ such that $\vert S_x\vert \leq 1$
for $x\in\R$, or if there is an $f$ with finite range that $S_x$
is always finite.

In this paper we solve some of those problems. We show that
there is always a uniformly antisymmetric $f:A\to\{0,1\}$ if
$A\subset \R$ is countable. We prove that the continuum hypothesis
 is equivalent to
the statement that there is an $f:\R\to\omega$ with $\vert
S_x\vert \leq 1$ for every $x\in\R$. If the continuum is at
least $\aleph_n$ then there  exists a point $x$ such that $S_x$
has at least $2^n-1$ elements. We also show that there is a
function $f:\Q\to\{0,1,2,3\}$ such that $S_x$ is always finite,
but no such function with finite range on $\R$ exists.
\medskip
\noindent{\bf Notation.} We use the standard set theory
notation. Notably, $\omega$ is the set of natural numbers,
ordinals are identified with the sets of smaller ordinals.
$\R$ is the set of reals, $\Q$ is the set of rationals. 
$\vert A\vert $ denotes the cardinality of $A$. If
$A$ is a set, $\kk$ is a cardinal, then $[A]^{\kk}=\{X\subseteq
A\colon \vert X\vert =\kk\}$, $[A]^{<\kk}=\{X\subseteq A\colon
\vert X\vert <\kk\}$. CH denotes the continuum hypothesis, i.e.,
that $\vert \R\vert =\aleph_1$.

\break
\medbreak
\centerline{\bf 1. Uniformly antisymmetric functions on countable sets}
\medskip
\proclaim Theorem 1. If $A\subseteq \R$ is countable, then there
is a uniformly antisymmetric function $f:A\to \{0,1\}$.

\P  Enumerate $A$ as $A=\{a_1,a_2,\dots\}$. By induction on
$n<\omega$ we define a finite set $\I_n=\{I_\gg:\gg\in\GG_n\}$
of open intervals such that
$\emptyset=\GG_0\subseteq\GG_1\subseteq\dots$, so 
$\emptyset=\I_0\subseteq \I_1\subseteq\dots$, each $I_\gg$ is of
the form 
$I_\gg=(b_\gg-h_\gg,b_\gg+h_\gg)$
with the following properties. Put $B_n=\{b_\gg:\gg\in\GG_n\}$.
\item{(1)} If $\gg\neq\gg'$ then either $I_\gg\cap
I_{\gg'}=\emptyset$, or one of them contains the other;
\item{(2)} if $I_{\gg'}\subseteq I_\gg$ then either
$I_{\gg'}\subseteq (b_\gg-h_\gg,b_\gg)$ or $I_{\gg'}\subseteq 
(b_\gg,b_\gg+h_\gg)$;
\item{(3)} $\{a_1,\dots,a_n\}\subseteq B_n$;
\item{(4)} $b_\gg\pm h_\gg\not\in A$ ($\gg\in\GG_n$) ;
\item{(5)} if we put $\varphi_\gg(x)=2b_\gg-x$ ($x\in I_\gg$,
$x\neq b_\gg$),
 then for
$I_{\gg'}\subseteq I_\gg$, $\varphi_\gg(I_{\gg'})\in \I_n$ holds.

To start, we put $\GG_0=\emptyset$.

If $\GG_{n-1}$ is already given, and $a_n\in B_{n-1}$, put
$\GG_n=\GG_{n-1}$. Otherwise, let $I_\gg$ be the unique shortest
interval in $\I_{n-1}$ containing $a_n$ if there exists one. 
Select
$I=(a_n-h,a_n+h)$ in such a way that it is either in
$(b_\gg-h_\gg,b_\gg)$ or in $(b_\gg,b_\gg+h_\gg)$ and
$\phi_{\gg_1}\cdots\phi_{\gg_r}(a_n\pm h)\not\in A$ for any
(applicable) product ($\gg_i\in \GG_{n-1}$). Notice that the
number of those products is $2^t$ where $t$ is the number of
intervals in $\I_{n-1}$ containing $a_n$. Now 
 add all $\phi_{\gg_1}\cdots\phi_{\gg_r}(I)$
to $\I_{n-1}$ and get $\I_n$. If no interval of $\I_{n-1}$
contains $a_n$ then let $I=(a_n-h,a_n+h)$, $a_n\pm h\not\in A$
be an arbitrary interval disjoint from those in $\I_{n-1}$ and
add it to get $\I_n$.

To conclude the proof of the Theorem we are going to show that
there exists a function $f:\R\to \{0,1\}$ such that
$f(\phi_\gg(x)) =1-f(x)$ ($\gg\in\bigcup\GG_n$). As $\phi^2_\gg$
is always a partial identity it suffices to show that no
$x\in\R$ is a fixed point of the product of odd many $\phi_\gg$.

Assume that $x=\phi_{\gg_1}\phi_{\gg_2}\cdots\phi_{\gg_t}(x)$, $t$
odd. Among the intervals $I_{\gg_1}, \dots,I_{\gg_t}$ there is a
longest one, say $I_\gg$ and that must contain all the others.
At every appearence of $\phi_\gg$ in the product 
$\phi_{\gg_1}\phi_{\gg_2}\cdots\phi_{\gg_t}$ the image of $x$ moves
from one side of $b_\gg$ to the other. $\phi_\gg$ therefore
appears even times. In the product the interval
$\phi_\gg\phi_{\gg_i}\cdots\phi_{\gg_j} \phi_\gg$ can be replaced by
$\phi_{\gg'_i}\cdots\phi_{\gg'_j}$ where
$I_{\gg'_r}=\phi_\gg(I_{\gg_r})$ ($i\leq r\leq j$), so
eventually we
succeed in eliminating an even number of $\phi$'s. We got a
shorter formula $x=\phi_{\gg'_1}\cdots\phi_{\gg'_{t'}}(x)$, but
$t'$ is still odd. Finally we get that $x=\phi^t_\gg(x)$ for
some odd $t$ which is impossible. \qed
\medskip
\medbreak
\centerline{\bf 2. When $S_x$ is finite}
\medskip
\noindent{\bf Definition.} If $f:\R\to\omega$ is a function,
then for $x\in \R$, set $S_x=\{h>0:f(x-h)=f(x+h)\}$.
\proclaim Theorem 2. There is a function $F\colon
[\omega_1]^{<\omega}\to \omega$ such that 
\item{(a)} if $F(A)=F(B)$ then
$\vert A\vert=\vert B\vert$;
\item{(b)} if $F(A)=F(B)$ then $A\cap B$ is an initial segment in
$A$, $B$; and 
\item{(c)} there do not exist $A_0$, $B_0$, $A_1$, 
$B_1\in [\omega_1]^{<\omega}$ such that $A_0\cup B_0=A_1\cup
B_1$, $F(A_0)=F(B_0)$, $F(A_1)=F(B_1)$, $A_0\neq B_0$, 
$A_1\neq B_1$, and $\{A_0,B_0\}\neq \{A_1,B_1\}$.

\P Let the diadic intervals of $\R $ be
$I_0,I_1,\dots$. For $\aa<\omega_1$ enumerate $\aa$ as
$\aa=\{\gg(\aa,i)\colon i<\omega\}$. (Recall that by our
axiomatic set theory assumptions $\aa$ is identified with the
set of smaller ordinals.) Select different irrational
numbers $r_{\aa}$ for $\aa<\omega_1$. We define a function
$c\colon [\omega_1]^2\to \omega$ as follows. We construct
$c(\bb,\aa)$ by induction on $\bb$, in the order of the
enumeration of $\aa$. For $\bb<\aa$, if $\bb=\gg(\aa,i)$, let
$c(\bb,\aa)$ be some $j<\omega$ such that

\item{(1)} $j>c(\gg(\aa,0),\aa),\dots,c(\gg(\aa,i-1),\aa)$ ;
\item{(2)} $r_{\bb}\in I_j$ ;
\item{(3)} $r_{\aa}\not\in I_j$ ;
\item{(4)} $r_{\xi}\not\in I_j$ for $\xi=\gg(\aa,0),\dots,\gg(\aa,i-1)$.

Clearly, such a $j<\omega$ can be found. Let, for $A\in
[\omega_1]^{<\omega}$, $F(A)$ be the isomorphism type of the
structure $(A;<,c)$, i.e., $F(A)=F(B)$ iff $\vert A\vert =\vert
B\vert$ and $c(a_i,a_j)=c(b_i,b_j)$ whenever $a_1<\cdots<a_n$,
$b_1< \cdots <b_n$ are the monotonic enumerations of $A$, $B$, respectively.
\medskip
\proclaim Claim 1. If $F(A)=F(B)$, then $A\cap B$ is an initial
segment in both sets.

\P Again, let $A=a_1,\dots,a_n$,
$B=b_1,\dots,b_n$ be the increasing enumerations. 
Assume that $a_i=b_j$ is a common element.
If $i\neq j$, say $i<j$, then $k=c(a_i,a_j)=c(b_i,b_j)$ has
$r_{a_i}\in I_k$ (by (2)), and $r_{b_j}\not\in I_k$ (by (3)), a
contradiction. So we have that $i=j$. If $t<i$, then, as
$c(a_t,a_i)=c(b_t,b_i)=c(b_t,a_i)$, $a_t=b_t$ by property (1). \qed
\medskip 
\proclaim Claim 2. There do not exist $\bb, \bb', \aa,
\aa'<\omega _1$
such that $\max(\bb,\bb')<\min(\aa,\aa')$, $c(\bb,\aa)=c(\bb',\aa')$, 
and $c(\bb',\aa)=c(\bb,\aa')$.

\P Set $i=c(\bb,\aa)$, $j=c(\bb',\aa)$. As $\bb,\bb'<\aa$,
$i\neq j$, say, $i<j$. Then, considering $c(\bb',\aa)$ we get
(by (4)) 
$r_{\bb}\not\in I_j$ while considering $c(\bb,\aa')$ we get that
$r_{\bb}\in I_j$, a contradiction. If $i>j$ we argue similarly. \qed

\medskip
Assume now that $F(A)=F(B)$ and we know $A\cup B$. We try to
reconstruct $A$, $B$. Put $X=A\cap B$, $Y=A-X$, $Z=B-X$. We can
assume that $m'=\max(Y)<\max(Z)=m$. In general, to every $x\in
Z$ let $x'$ be the element in $Y$ corresponding to $x$ under the
(unique) order-preserving bijection between $Z$ and $Y$. 

For $a<b$ in $A$, $c(a,b)$ codes a diadic interval including
$r_a$ but excluding $r_b$. The structure $(A;<,c)$ gives a
diadic interval for every element in $A$ separating it from the
rest of $A$. As $F(A)=F(B)$ this interval is the same for $x$
and $x'$. We get therefore, that there is a diadic interval
containing $r_x$, $r_{x'}$ but nothing else from $A\cup B$. This
makes possible to find $x'$ if $x$ is given, or to find $x$ if
$x'$ is given. Anyway, we can find $m'$.

\medskip
\proclaim Claim 3. $X=\{x\in A\cup B\colon x<m'$ and $c(x,m')=c(x,m)\}$.

\medskip
\P $\subseteq$ is obvious. If, say $x\in Z$ and $c(x,m')=c(x,m)$
then $c(x,m')=c(x,m)=c(x',m')$ a contradiction to (1), as $x\neq
x'$. \qed

\medskip
As now $X$ is known, we can decompose $Y\cup Z$ into the 
pairs $\{x,x'\}$ by the
argument before Claim 3. Given such a pair $\{u,v\}$ we have to
find if $u\in Z$, $v\in Y$ or vice versa. We know that
$c(x',m')=c(x,m)$, so, knowing $m$, $m'$ we can identify $x$,
$x'$ if we can show that $c(x,m')\neq c(x',m)$. But this is done
in Claim 2. \qed

\medskip
\proclaim Theorem 3. If CH holds, then there is a function
$f\colon \R \to \omega$ such that for every $x\in \R$ 
$
S_x$ has at most one element.

\P Let $\{b_{\aa}\colon \aa<\omega_1\}$ be a Hamel basis,
$F:[\omega]^{<\omega}\to \omega$ as
in Theorem 1. To 
$$
x=\sum^n_{i=1}\ll_i b_{\aa_i}
$$
($\ll_i\neq 0$, $\ll_i\in\Q$), $\aa_1<\cdots<\aa_n$ we associate
some $f(x)$ that codes the ordered string $\langle
\ll_1,\dots,\ll_n\rangle$ as well as $F(\{\aa_1,\dots,\aa_n\})$ .
This is possible as there are countably many possibilities for
both. 

Assume that $x\neq y$, $f(x)=f(y)$. We try to recover the pair
$\{x,y\}$ from $x+y$. By our coding of the string of the
coefficients in the Hamel basis and the properties of the function $F$
described in the previous Theorem, $x$, $y$ can be written as
$$
x=\sum^n_{i=1}\ll_ib_{\aa_i},\quad y=\sum^n_{i=1}\ll_ib_{\bb_i}
$$
such that $\aa_i=\bb_i$ for $1\leq i\leq m$ (some $m< n$),
and
$\{\aa_{m+1},\dots,\aa_n\}\cap\{\bb_{m+1},\dots,\bb_n\}=\emptyset$.
$x+y$ can be written in the above basis as
$$
x+y=\sum^m_{i=1}(2\ll_i)b_{\aa_i} +\sum^n_{i=m+1}\ll_ib_{\aa_i}+
\sum^n_{i=m+1}\ll_ib_{\bb_i}.
$$
The support of $x+y$, i.e., the set of those basis vectors in
which it has nonzero coefficients is
$$
\{\aa_1,\dots,\aa_m,\aa_{m+1},\dots,\aa_n,\bb_{m+1},\dots,\bb_n\}
$$
from which, by the previous Theorem $\{\aa_1,\dots,\aa_n\}$ and
$\{\bb_1,\dots,\bb_n\}$ can be recovered. Then we can find
$\ll_1,\dots,\ll_n$, i.e., $x$ and $y$ can be reconstructed. \qed

\medskip
Before proving that if a vector space $V$ with $\vert V\vert
\geq \omega_n$ is $\omega$-colored then $\vert S_x\vert\geq
2^n-1$ holds 
for some $x\in V$ we give a proof of the combinatorial part of
the theorem. We then show how to modify the proof to get the
stated result.
\medskip
\proclaim Theorem 4. If $2\leq n<\omega$ and 
$f:[\omega_n]^{<\omega}\to \omega$ then there exists a set $s\in
[\omega_n]^{<\omega}$ which can be written in $2^n-1$ ways as the
union of two different sets $s=x\cup y$ such that $f(x)=f(y)$.

\P Assume that $f:[\omega_n]^{<\omega}\to \omega$. Select
$\omega_{n-1}< y^0_n<\omega_n$ such that it is not in any of the
sets
$$
\{x:\omega_{n-1}<x<\omega_n, f(s_1\cup\{x\})=j_1,\dots,f(s_t\cup\{x\})
=j_t\}
$$
(for some $s_1,\dots,s_t\in [\omega_{n-1}]^{<\omega}$,
$j_1,\dots,j_t<\omega$) which happen to have  one
element. This is possible, as the number of those sets is
$\aleph_{n-1}$, and they are all small enough. 

Assume now that $y^0_{i+1},\dots,y^0_n$ are already defined. Let
$\omega_{i-1}<y^0_i<\omega_i$ be such that it is not in any of
the sets of the form 
$$
\{x: f(s_1\cup\{x,y^0_{i+1},\dots,y^0_n\})
=j_1,\dots,f(s_t\cup \{x,y^0_{i+1},
\dots,y^0_n\})=j_t, \omega_{i-1}<x<\omega_i\}
$$
for some $s_1,\dots,s_t\in [\omega_{i-1}]^{<\omega}$,
$j_1,\dots,j_t<\omega$, which are singletons. Again,
this choice is possible. 

If $y^0_1,\dots,y^0_n$ are given, we define $y^1_i$ ($1\leq
i\leq n$) in increasing order. Select $y^1_1\neq y^0_1$ such that
$\omega<y^1_1<\omega_1$ and
$f(\{y^1_1,y^0_2,\dots,y^0_n\})=f(\{y^0_1,\dots,y^0_n\})$. This
is possible, as otherwise $y^0_1$ would be the only element in 
$\{x:\omega< x<\omega_1, f(\{x,y^0_2,\dots,y^0_n\})=j\}$ where 
$j=f(\{y^0_1,y^0_2,\dots,y^0_n\})$, a contradiction to the
choice of $y^0_1$. 

If $y^1_1,\dots,y^1_{i-1}$ have already been selected, let
$y^1_i\neq y^0_i$ be such that $\omega_{i-1}<y^1_i<\omega_i$ and 
$$
f(s\cup \{y^1_i,y^0_{i+1},\dots,y^0_n\})=f(s\cup 
\{y^0_i,y^0_{i+1},\dots,y^0_n\})
$$
for every $s\subseteq
\{y^0_1,y^1_1,\dots,y^0_{i-1},y^1_{i-1}\}$. This is possible by
the choice of $y^0_i$. 
\medskip
For $1\leq k\leq m\leq n$, $g:\{k,\dots,m\}\to \{0,1\}$ put 
$A=\{y^0_1,y^1_1,\dots,y^0_{k-1},y^1_{k-1}\}$,
$B=\{y^0_k,\dots,y^0_m\}$,
$B^g=\{y^{g(k)}_k,\dots,y^{g(m)}_m\}$,
$C=\{y^0_{m+1},\dots,y^0_n\}$. 
\medskip
\proclaim Claim. $f(A\cup B\cup C)=f(A\cup B^g\cup C)$.

\P By induction on $m$. The inductive step trivially follows
from the choice of $y^1_m$. \qed
\medskip
To conclude the proof of the Theorem, assume that $1\leq k\leq
n$, $g:\{k,\dots,n\}\to \{0,1\}$. Put 
$A=\{y^0_1,y^1_1,\dots,y^0_{k-1},y^1_{k-1}\}$,
$B^g=\{y^{g(k)}_k,\dots,y^{g(n)}_n\}$ and let $1-g$ be the
function with $(1-g)(i)=1-g(i)$ for $k\leq i\leq n$. Using the
Claim we get that $f(A\cup B^g)=f(A\cup B^{1-g})$ and clearly 
$(A\cup B^g)\cup (A\cup
B^{1-g})=\{y^0_1,y^1_1,\dots,y^0_n,y^1_n\}$. The number of those
decompositions, i.e., that of the pairs $\{g,1-g\}$ is
$2^{n-k}$, summing we get  $2^{n-1}+\cdots+1=2^n-1$. \qed

\medbreak
\proclaim Theorem 5. Let $V$ be a vector space, $\vert V\vert
\geq \aleph_n$ $(2\leq n<\omega)$ and $f:V\to \omega$ be given.
Then $\vert S_x \vert \geq 2^n-1$ for some $x\in V$.

\P Assume that $\{b(\aa):\aa<\omega_n\}$ is a linearly
independent set. Select $\omega_{n-1}<y^0_n<\omega_n$ outside
any of the one-element sets of the form
$$\displaylines{
\Bigl\{
\omega_{n-1}<x<\omega_n:f\Bigl( \sum_{z\in s_1}b(z)+
{1\over 2}\sum_{z\in s_1'}b(z)+b(x)
\Bigr)=j_1,\dots,\hfill\cr
\hfill f\Bigl( \sum_{z\in s_t}b(z)+
{1\over 2}\sum_{z\in s_t'}b(z)+b(x)
\Bigr)=j_t
\Bigr\}\cr}
$$
where $s_1,s_1',\dots,s_t,s_t'\in [\omega_{n-1}]^{<\omega}$, $j_1\dots,j_t<
\omega$. Given $y^0_{i+1},\dots,y^0_n$ , let $\omega_{i-1}<y^0_i<
\omega_i$ be not in any of the one-element sets
$$\displaylines{
\Bigl\{
\omega_{i-1}<x<\omega_i:f\Bigl( \sum_{z\in s_1}b(z)+
{1\over 2}\sum_{z\in s_1'}b(z)+b(x)+b(y^0_{i+1})+\cdots+b(y^0_n)
\Bigr)=j_1,\dots,\hfill\cr
\hfill f\Bigl( \sum_{z\in s_t}b(z)+
{1\over 2}\sum_{z\in s_t'}b(z)+b(x)+b(y^0_{i+1})+\cdots+b(y^0_n)
\Bigr)=j_t
\Bigr\}\cr}
$$
where $s_1,s_1',\dots,s_t,s_t'\in [\omega_{i-1}]^{<\omega}$, $j_1,\dots,j_t<
\omega$. If $y^0_1,\dots,y^0_n$ are already constructed, let
$y^1_1\neq y^0_1$ be such that $\omega<y^1_1<\omega_1$ and 
$f\bigl(b(y^1_1)+b(y^0_2)+\cdots+b(y^0_n)\bigr)=
f\bigl(b(y^0_1)+b(y^0_2)+\cdots+b(y^0_n)\bigr)$. With
$y^1_1,\dots,y^1_{i-1}$  defined, let
$\omega_{i-1}<y^1_i<\omega_i$, $y^1_i\neq y^0_i$ be such that for
every $s\cup s'\subseteq
\{y^0_1,y^1_1,\dots,y^0_{i-1},y^1_{i-1}\}$, if $s\cap
s'=\emptyset$, then 
$$
\displaylines{
f\Bigl( \sum_{z\in s}b(z)+
{1\over 2}\sum_{z\in s'}b(z)+b(y^1_i)+b(y^0_{i+1})+\cdots+b(y^0_n)
\Bigr)=\hfill\cr
\hfill 
f\Bigl( \sum_{z\in s}b(z)+
{1\over 2}\sum_{z\in s'}b(z)+b(y^0_i)+b(y^0_{i+1})+\cdots+b(y^0_n)
\Bigr)
\cr}
$$
holds. This is possible by the choice of $y^0_i$.

For $1\leq k\leq m\leq n$, $g:\{k,\dots,m\}\to \{0,1\}$ we define
$$\displaylines{
A_k={1\over
2}\bigl(b(y^0_1)+b(y^1_1)+\cdots+b(y^0_{k-1})+b(y^1_{k-1})\bigr),\cr
B=b(y^0_k)+\cdots+b(y^0_m),\quad B^g
=b(y^{g(k)}_k)+\cdots+b(y^{g(m)}_m),\cr
C=b(y^0_{m+1})+\cdots+b(y^0_n).\cr}
$$
\medskip
\proclaim Claim. $f(A_k+B+C)=f(A_k+B^g+C)$.

\P As in Theorem 4. \qed
\medbreak
To conclude the proof one can argue as in Theorem 4, and
decompose $b(y^0_1)+b(y^1_1)+\cdots+b(y^0_n)+b(y^1_n)$  in $2^n-1$
ways into the sum of two vectors with the same $f$ value as 
$(A_k+B^g)+(A_k+B^{1-g})$ where $g:\{k,\dots,n\}\to\{0,1\}$. \qed
\medbreak
\centerline{\bf 3. Finite range}
\medskip
\proclaim Theorem 6. There is a function $f:\Q\to \{0,1,2,3\}$
such that for every $x\in\Q$, $S_x$ is finite.

\P It suffices to find such a function assuming  two values on 
the set 
$\Q^+=\{x\in\Q:x>0\}$. We decompose $\Q^+$ into the increasing
union of finite sets $A_1\subseteq A_2\subseteq\cdots$. We also
define an auxiliary graph $G$ on $\Q^+$. Two points $x$ and $y$
will be joined in $G$ if $(x+y)/2\in A_n$ for some $n$ but $x$,
$y\not\in A_{n+1}$. If, with an appropriate choice of the sets
we can guarantee that the graph $G$ is bipartite, then the
bipartition of $G$ will give a good function on $\Q^+$. Indeed,
if $x\in A_n$ and $f(x-h)=f(x+h)$ then one of $x-h$, $x+h$ is in
$A_{n+1}$ so there are only finitely many such $h$'s. 

Let a positive rational number be in $A_n$ if it is of the form
$x=p/n!$ and $x<2^n$. Clearly these sets are finite, constitute
an increasing sequence, and their union is $\Q^+$.

We first show that if $x$, $y$ are joined in $G$, then they
first appear in the same $A_n$. Assume that $x\in A_{n+1}-A_n$, $y\in
A_{m+1}-A_m$, $m\geq n$, and $z=(x+y)/2\in A_{n-1}$. Then, the
denominator of $y=2z-x$ is (a divisor of) $(n+1)!$. Also,
$y<2z<2^{n+1}$, so $y\in A_{n+1}$, i.e., $m=n$.

Finally, we  show that $G$ on $A_{n+1}-A_n$ does not
contain odd circuits. Assume that $a_1,\dots,a_{2u+1}$ is one,
i.e., $a_i+a_{i+1}=2b_i$ for some $b_i\in A_{n-1}$ ($1\leq i
\leq 2u+1$). Here, we use cyclical indexing, i.e.,
$a_{2u+2}=a_1$. Then again, $a_1<2b_1<2^n$, and as 
$a_1=b_1-b_2+b_3-\cdots+b_{2u+1}$, $a_1$ has denominator
$(n-1)!$, so it is in $A_n$, a contradiction. \qed

\proclaim Theorem 7. If $f:\R\to \{1,\dots,n\}$ is a function, then $S_x$
is infinite for some $x\in \R$.

\P Actually the result is true for any uncountable vector space
$V$ 
over $\Q$. Assume that $f:V\to \{1,\dots,n\}$. 
Let $\{b(\aa):\aa<\omega_1\}$ be linearly independent. For
$\bb<\aa<\omega_1$, 
the formula $F(\bb,\aa)=f\bigl(b(\aa)-b(\bb)\bigr)$ defines an
$n$-coloring of $[\omega_1]^2$. By an old Erd\H os-Rado theorem
(see Cor.~1, p.~459 in [3]),
there are a color $1\leq k\leq n$ and ordinals 
$\aa(0)<\cdots<\aa(\omega)$, such that $F\bigl(\aa(i),\aa(j)\bigr)=k$ for
$i<j\leq \omega$. But then, 
$$
f\bigl(b(\aa(i))-b(\aa(0))\bigr)=
f\bigl(b(\aa(\omega)\bigr)-b\bigl(\aa(i))\bigr)=k,
$$
i.e., the vector $b\left(\aa(\omega)\right)-b\left(\aa(0)\right)$ 
can be written
infinitely many ways as the sum of two monocolored vectors. \qed

\medbreak
\centerline{\bf References}
\medskip
\item{[1]} K. Ciesielski: Notes on  problem 1 from ``Uniformly antisymmetric
functions'', to appear.
\item{[2]} K. Ciesielski, L. Larson: Uniformly antisymmetric
functions, to appear.
\item{[3]} P. Erd\H os, R. Rado: A partition calculus in set
theory, {\sl Bull. of the Amer. Math. Soc. \bf 62} (1956), 427--489.
\item{[4]} P. Komj\'ath: Vector sets with no repeated differences, {\sl Coll.
Math. \bf 64}(1993), 129--134.
\item{[5]} P. Kostyrko: There is no strongly locally antisymmetric
set,  {\sl Real Analysis Exchange, \bf 17} (1991/92), 423--425.
\item{[6]} B.~S.~Thomson: {\sl Symmetric properties of real
functions}, to appear.

\bye